\documentclass[12pt]{article}
\usepackage{amsmath,amsthm, amssymb}
\def\Omega{\varOmega}
\def\ep{\varepsilon}

\begin{document}

\title{\large 
\bf A centennial of the Zaremba--Hopf--Oleinik Lemma}
\author { {\it A.I.~Nazarov}\footnote {Partially supported
by RFBR grant 09-01-00729 and by grant NSh.4210.2010.1.},\\ 
Saint-Petersburg University, \\
{\small e-mail: \ al.il.nazarov@gmail.com}} 
\date{}\maketitle

{\small We consider the Hopf--Oleinik normal derivative lemma for elliptic and parabolic equations 
under minimal restrictions on lower-order coefficients. Boundary gradient estimates for solutions
are also established.}

\section{Introduction}

Qualitative theory of partial differential equations is in intensive 
development over last half of century. In this paper we discuss the Hopf--Oleinik Lemma, 
one of the most important tools in studying solutions to elliptic and parabolic equations, in particular,
the key argument in the proof of uniqueness theorems.\medskip

For the Laplace operator this property is well known 
for one hundred years, starting from a pioneer paper of Zaremba \cite{Z}, and reads as follows. Let $\partial\Omega\in{\cal C}^2$ and
let ${\cal L}=-\Delta$. Then, if $0\in\partial \Omega$, we have
$${\cal L}u=f\ge0\ \ \mbox{in}\ \ \Omega;\quad
u(x)>u(0)\ \ \mbox{in}\ \  \Omega\quad \Longrightarrow\quad
\frac {\partial u}{\partial\bf n}(0)<0.\eqno({\bf ZHO})$$

For general operators of non-divergence type with bounded measurable coefficients this result was established 
in elliptic case independently by E.~Hopf \cite{Ho} and O.A.~Oleinik \cite{O} and in parabolic case 
by L.~Nirenberg \cite{Ni}. Later the efforts of many mathematicians were aimed at the reduction of the boundary 
smoothness\footnote{See also an earlier paper \cite{G} for equations with H\"older continuous leading coefficients.}. 
They established that the sharp condition for ({\bf ZHO})
to fulfil is the Dini condition for the boundary normal, see, e.g., \cite{Hi}. In a weakened
form (the existence a boundary point $x^0$ in any neighborhood of the origin and a direction $\ell$ 
such that $\frac {\partial u}{\partial\ell}(x^0)<0$) this fact holds true for a much wider class
of domains including all Lipschitz ones, see \cite{Na} for elliptic 
equations and \cite{K} for parabolic ones. Note that all these results are related to classical 
solutions, i.e. $u\in {\cal C}^2(\Omega)$ in elliptic case and $u\in {\cal C}^{2,1}(Q)$ in 
parabolic case.\medskip

Now let us consider generalized ({\it strong}) solutions for non-divergence type equations
$${\cal L}u\equiv-a_{ij}(x)D_iD_ju+b_i(x)D_iu=f(x);\eqno({\bf NDE})$$
$${\cal M}u\equiv \partial_tu-a_{ij}(x;t)D_iD_ju+b_i(x;t)D_iu=f(x;t),\eqno({\bf NDP})$$
i.e. we assume $D(Du)\in L_{n,loc}(\Omega)$ in ({\bf NDE}) and
$\partial_tu$, $D(Du)\in L_{n+1,loc}(Q)$ in ({\bf NDP}) (in the parabolic case
also some anisotropic spaces are admissible).\medskip

We always suppose that operators under consideration are uniformly elliptic 
(parabolic), i.e. for all values of arguments 
\begin{equation}\label{ell}
\nu|\xi|^2\le a_{ij}(\cdot)\xi_i\xi_j\le\nu^{-1}|\xi|^2,\qquad 
\xi\in{\mathbb R}^n,
\end{equation}
where $\nu$ is a positive constant. Note that we can also assume $a_{ij}\equiv a_{ji}$ without loss of generality.\medskip

The properties of generalized solutions to the equations ({\bf NDE})-({\bf NDP}), 
under assumption that leading coefficients $a_{ij}$ are only measurable,
were investigated in a number of papers\footnote{We mention in this connection a quite recent paper \cite{M11} discussing some degenerate elliptic equations.}. The problem of our interest is how 
``bad'' may be lower-order coefficients $b_i$ to ensure the Hopf--Oleinik Lemma 
to hold true. We provide sharp conditions for this. We also touch the topic closely related 
to ({\bf ZHO}), especially in idea of proof, 
namely, the gradient estimates at the boundary.\medskip

Note that for divergence type equations 
$$-D_i\big(a_{ij}(x)D_ju\big)+b_i(x)D_iu=0;\eqno({\bf DE})
$$
$$\partial_tu-D_i\big(a_{ij}(x;t)D_ju\big)+b_i(x;t)D_iu=0.\eqno({\bf DP})
$$
({\bf ZHO}) does not hold. 
The simplest counterexample is the function $u(x)=x_2^2+2x_2|x_1|$ which
is positive in the upper half-plane, satisfies the equation ({\bf DE}) with
$$(a_{ij})=\bmatrix 1&-\mbox{sign}(x_1)\\-\mbox{sign}(x_1)&2\endbmatrix;\qquad
b_i\equiv0;\qquad f\equiv0,
$$
but $u(0,0)=0$ and $D_2u(0,0)=0$.

Moreover, even continuity of $a_{ij}$ does not improve the situation.
Let us describe corresponding counterexample\footnote{See also \cite[Problem 3.9]{GT}.}.

Let $\Omega$ be a convex domain, and let $0\in\partial\Omega$. Assume that at 
the neighborhood of the origin $\partial\Omega$ is the graph of a function 
$x_n=\phi(x')$. Finally, suppose that $\phi\in{\cal C}^1$ but $D'\phi$ 
is not Dini continuous at the origin. 

As it was mentioned the Hopf--Oleinik lemma for the Laplacian fails 
in such domain. Now we rectify the boundary near the origin and obtain an 
operator of the form ({\bf DE}) with {\it continuous} leading coefficients 
and $b_i\equiv0$ for which ({\bf ZHO}) fails in smooth domain. Considering
functions depending only on spatial variables we see that this example works 
also for the parabolic operator ({\bf DP}).\medskip

The paper is organized as follows. In Section 2 we deal with elliptic equations,
Section 3 is devoted to parabolic equations. In both sections we use the ``composite''
variant of the A.D.~Aleksandrov maximum estimate (\cite{Li00}; see also \cite{AN} for a weaker version)
and slightly modify classical techniques due to Ladyzhenskaya--Ural'tseva \cite{LU88}, see also \cite{S10}.\medskip 


Let us recall some notation. $x=(x_1,\dots,x_{n-1},x_n)=(x',x_n)$ is a vector 
in $\mathbb R^n$, $n\ge2$, with the Euclidean norm $|x|$; $(x;t)$ is a point 
in $\mathbb R^{n+1}$. 

$\Omega$ is a domain in $\mathbb R^n$ and 
$\partial\Omega$ is its boundary; ${\text{\bf n}}=({\text{\bf n}}_i(x))$ is 
the unit vector of the outward normal to $\partial\Omega$ at the point $x$.

For a cylinder $Q=\Omega\times]0,T[$ we denote by
$\partial''Q=\partial\Omega\times]0,T[$ its lateral surface and by
$\partial'Q=\partial''Q\cup\{\overline{\Omega} \times \{0\}\}$ its
parabolic boundary.


We define
$$\begin{array}{ll}
B_r(x^0)=\{x\in \mathbb R^n : |x-x^0|<r\}, & B_r=B_r(0);\\
{\cal B}_{r,h}({x^0}\vphantom{x}')=\{x\in \mathbb R^n: |x'-{x^0}\vphantom{x}'|<r, 0<x_n<h\}; &
{\cal B}_{r,h}={\cal B}_{r,h}(0);\\
Q_r(x^0;t^0)=B_r(x^0)\times]t^0-r^2;t^0[, & Q_r=Q_r(0;0);\\
{\cal Q}_{r,h}({x^0}\vphantom{x}';t^0)={\cal B}_{r,h}({x^0}\vphantom{x}')\times]t^0-r^2;t^0[, & 
{\cal Q}_{r,h}={\cal Q}_{r,h}(0;0).
\end{array}
$$

The indices $i, j$  vary from 1 to $n$.
Repeated indices indicate summation. 

The symbol $D_i$ denotes the operator of differentiation with respect to 
$x_i$; in particular, $Du=(D_1u, \dots, D_nu)=(D'u,D_nu)$ is the gradient of 
$u$. $\partial_tu$ stands for the derivative of $u$ with respect to $t$.

We denote by $\|\cdot\|_{q,\Omega}$ the norm in $L_q(\Omega)$.
We introduce two scales of anisotropic spaces:

$L_{q,\ell}(Q)=L_{\ell}\big(\,]0,T[\,\to L_q(\Omega)\big)$ with the norm
$\|f\|_{q,\ell,Q}=\big\|\|f(\cdot;t)\|_{q,\Omega}\big\|_{\ell,]0,T[}$;

$\widetilde L_{q,\ell}(Q)=L_q\big(\Omega\to L_{\ell}(\,]0,T[\,)\big)$ with the norm
$\|f\|^{\sim}_{q,\ell,Q}=\left\|\|f(x;\cdot)\|_{\ell,]0,T[}\right\|_{q,\Omega}$.

\noindent Obviously, $L_{q,q}(Q)=\widetilde L_{q,q}(Q)=L_q(Q)$. Further, by the Minkowskii inequality,
$$\|f\|^{\sim}_{q,\ell,Q}\le\|f\|_{q,\ell,Q}\quad\mbox{for}\ \ q\ge \ell;\qquad
\|f\|_{q,\ell,Q}\le\|f\|^{\sim}_{q,\ell,Q}\quad\mbox{for}\ \ q\le \ell.
$$
We denote by $\widehat L_{q,\ell}(Q)$ the space
$$L_{q,\ell}(Q)\cap\widetilde L_{q,\ell}(Q)=
\begin{cases}L_{q,\ell}(Q),& q\ge \ell;\\
\widetilde L_{q,\ell}(Q),& q\le \ell
\end{cases}
$$
with the norm $|\!|\!|f|\!|\!|_{q,\ell,Q}=\max\{\|f\|_{q,\ell,Q},\|f\|^{\sim}_{q,\ell,Q}\}$.\medskip

{\bf Remark 1}. Note that we always deal with the space $\widehat L_{q,\ell}(Q)$ i.e. take the more strong of two norms. The reason is that up to now anisotropic versions of the Aleksandrov--Krylov maximum principle
 (see \cite{N87}, \cite{N01}) are proved only in terms of stronger norm.\medskip

We set $f_+=\max \{f,0\},\quad f_-=\max \{-f,0\}$.

Following \cite{Li00}, we say that $\omega:\ [0,1]\to\mathbb R_+$ belongs to the class ${\cal D}_1$ if
$\omega(1)=1$, $\omega$ is continuous and increasing while $\omega(\sigma)/\sigma$
is summable and decreasing. In this case we define 
${\cal I}_{\omega}(s)=\int_0^s\frac {\omega(\sigma)}{\sigma}\,d\sigma$.

We use letters $M$, $N$, $C$ (with or without indices) to denote various 
constants. To indicate that, say, $N$ depends on some parameters, we list 
them in the parentheses: $N(\dots)$.

\section{Elliptic case}

Recall that in this section we assume $D(Du)\in L_{n,loc}(\Omega)$. 

The next statement is a particular case of \cite[Theorem 3.2]{Li00}.\medskip

{\bf Proposition 2.1}.  {\it Let $\cal L$ be an operator of the form (\,{\bf NDE}) in a bounded, strictly Lipschitz domain $\Omega$, and let the condition (\ref{ell}) be satisfied. Suppose also that the vector function 
${\bf b}=(b_i)$ can be written as follows:
\begin{equation}\label{decomp1}
{\bf b}={\bf b}^{(1)}+{\bf b}^{(2)};\qquad |{\bf b}^{(1)}|\in L_n(\Omega);
\end{equation}
\begin{equation*}
|{\bf b}^{(2)}|\le {\mathfrak B}\,\frac {\omega(d/{\rm diam}(\Omega))}d,\qquad\omega\in {\cal D}_1,
\end{equation*}
where $d=d(x)={\rm dist}(x,\partial\Omega)$.

Then for any solution of ${\cal L}u=f
$ in $\Omega$ satisfying $u\big|_{\partial\Omega}\le 0$, 
the following estimate holds:
$$u\le N_0\cdot\frac{{\rm diam}(\Omega)}{\nu}\cdot \|f
_+\|_{n,\{u>0\}}
,$$
provided $\|{\bf b}^{(1)}\|_{n,\Omega}\le{\mathfrak B}_0$, where $N_0$ and ${\mathfrak B}_0$ depend only on 
$n$, $\nu$, ${\mathfrak B}$ and the Lipschitz constant of $\partial\Omega$.}
\medskip

Now we prove a quantitative version of tne maximum principle, the so-called ``boundary growth lemma'' 
(its versions for $|{\bf b}|\in L_n$ are proved, e.g., in \cite[Lemma 2.5']{LU85} and 
\cite[Lemma 2.6]{S10}).\medskip

{\bf Lemma 2.2}. {\it Let $\cal L$ be an operator of the form (\,{\bf NDE}) in 
${\cal B}_{\rho,\rho}$, $\rho\le R$, and let the conditions (\ref{ell}), (\ref{decomp1}) and
\begin{equation}\label{decomp2}
|{\bf b}^{(2)}(x)|\le {\mathfrak B}\,\frac {\omega(x_n/R)}{x_n},\qquad  \omega\in {\cal D}_1,
\end{equation}
 be satisfied. Suppose also that $\|{\bf b}^{(1)}\|_{n,{\cal B}_{\rho,\rho}}\le {\mathfrak B}_0$ where 
${\mathfrak B}_0={\mathfrak B}_0(n,\nu, {\mathfrak B})$ is the constant from Proposition 2.1. 
If $u$ is a nonnegative solution of ${\cal L}u=f\ge0$ in ${\cal B}_{\rho,\rho}$ satisfying $u\ge k$ 
on $\partial {\cal B}_{\rho,\rho}\cap\{x_n=0\}$ for some $k>0$ then for $\xi\le\frac 12$ the inequality}
\begin{equation}\label{lem22}
u\ge k\cdot\big(\beta(n,\nu,{\mathfrak B},\omega,\xi)-
N_1(n,\nu,{\mathfrak B},\omega,\xi)\cdot(\|{\bf b}^{(1)}\|_{n,{\cal B}_{\rho,\rho}}+{\mathfrak B}\omega(\rho/R))\big)
\end{equation}
{\it holds in ${\cal B}_{(1-\xi)\rho,(1-\xi)\rho}$ with some positive constants $\beta$ and $N_1$.
}\medskip

{\bf Proof}. Consider the barrier function
\begin{equation}\label{barrier}
w(x)=\Big(1-A\,\frac{x_n}{\rho}\Big)^2
+2(1+A)\,\big(\varphi(x_n/\rho)-\varphi(1/A)\big)-\frac {|x'|^2}{\rho^2},
\end{equation}
where (cf. \cite{Li00})
$$\varphi(s)=
\int\limits_0^s\Big(\exp\Big(\frac{\mathfrak B}{\nu}\,{\cal I}_{\omega}(\sigma)\Big)-1\Big) d\sigma,
$$
while $A\ge\frac {\sqrt{n-1}}{\nu}$ is a constant to be defined later.\medskip

Direct calculation shows that for $x\in{\cal B}_{\rho,\frac{\rho}A}$
$$-a_{ij}D_iD_jw\le\frac 2{\rho^2}\cdot\big[-\nu A^2-\nu(1+A)\varphi''(\tau)+(n-1)\nu^{-1}\big]
\le -\frac {2\nu}{\rho^2}\,(1+A)\varphi''(\tau);
$$
$$b^{(1)}_iD_iw\le |{\bf b}^{(1)}|\cdot\frac 2{\rho}\cdot\big[A\big(1-A\,\tau\big)
+(1+A) \varphi'(\tau)+1\big]\le |{\bf b}^{(1)}|\cdot\frac 2{\rho}\,(1+A)(1+\varphi'(\tau));
$$
$$b^{(2)}_iD_iw\le\frac 2{\rho^2}\cdot\big[A\big(1-A\,\tau\big)
+(1+A) \varphi'(\tau)+1\big]\cdot{\mathfrak B}\,\frac {\omega(\tau)}{\tau}\le
\frac {2\mathfrak B}{\rho^2}\,(1+A)\,\frac {\omega(\tau)}{\tau}\,(1+\varphi'(\tau))
$$
(here $\tau=x_n/R$). 
Since $\varphi''(\tau)=\frac {\mathfrak B}{\nu}\,\frac {\omega(\tau)}{\tau}\,(1+\varphi'(\tau))$, we have
${\cal L}w\le C_1(A,\nu,{\mathfrak B})|{\bf b}^{(1)}|\rho^{-1}$.\medskip

Further, $w(x)\le0$
for $|x'|=\rho$, $0<x_n<\frac{\rho}A$ and for $|x'|\le\rho$, $x_n=\frac{\rho}A$. Finally, $w(x)\le1$
for $|x'|\le\rho$, $x_n=0$. This gives $kw-u\le0$ on $\partial{\cal B}_{\rho,\frac{\rho}A}$.\medskip

Proposition 2.1 gives for $x\in{\cal B}_{\rho,\frac{\rho}A}$
$$u(x)\ge kw(x)-kC_2\rho\cdot\|({\cal L}w)_+\|_{n,{\cal B}_{\rho,\frac{\rho}A}}\ge
kw(x)-kC_3\|({\bf b}^{(1)}\|_{n,{\cal B}_{\rho,\frac{\rho}A}},
$$
where $C_2$ depends only on $n$, $\nu$ and ${\mathfrak B}$ while $C_3$ depends on the same quantities and on
$A$.\medskip

Now we observe that $\varphi(s)=o(s)$ as $s\to0$. Thus, we can choose 
$A=A(n,\nu,{\mathfrak B},\omega,\xi)\ge\frac {\sqrt{n-1}}{\nu}$ so large that $2(1+A)\varphi(1/A)\le\frac {\xi}2$. Since 
$(1-\frac {\xi}2)^2-\frac {\xi}2-(1-\xi)^2=\xi(1-\frac {3\xi}4)\ge\frac {5\xi}8$, this gives 
\begin{equation}\label{q}
u\ge k\cdot\Big(\frac {5\xi}8-
C_4(n,\nu,{\mathfrak B},\omega,\xi)\|{\bf b}^{(1)}\|_{n,{\cal B}_{\rho,\rho}}\Big)_+
\equiv k_1\qquad\mbox{in}\quad {\cal B}_{(1-\xi)\rho,\frac {\xi\rho}{2A}}.
\end{equation}

Now we consider the set 
$K_{\rho}={\cal B}_{\rho,\rho}\setminus{\cal B}_{\rho,\frac {\xi\rho}{4A}}$. Note that 
coefficients $b^{(2)}_i$ are bounded on this set, and
\begin{equation}\label{qq}
\|{\bf b}^{(2)}\|_{n,K_{\rho}}\le C_5(n){\mathfrak B}\Big(\int\limits_{\frac {\xi}{4A}}^1
\Big(\frac {\omega(s\rho/R)}s\Big)^n ds\Big)^{\frac 1n}\le C_6(n,\nu,{\mathfrak B},\omega,\xi)
{\mathfrak B}\omega(\rho/R).
\end{equation}
We apply ``the ink-spot expansion lemma'' (\cite[Lemma 2.2]{LU85}) and obtain
$$u\ge k_1\cdot\Big(\varkappa(n,\nu,{\mathfrak B},\omega,\xi)-
C_7(n,\nu,{\mathfrak B},\omega,\xi)\|{\bf b}^{(1)}+{\bf b}^{(2)}\|_{n,K_{\rho}}\Big)\qquad\mbox{in}\quad 
{\cal B}_{(1-\xi)\rho,(1-\xi)\rho}\setminus{\cal B}_{(1-\xi)\rho,\frac {\xi\rho}{2A}}.
$$
By (\ref{q}) and (\ref{qq}) we arrive at (\ref{lem22}).\hfill$\square$\medskip

{\bf Lemma 2.2'}. {\it Let $\cal L$ be as in Lemma 2.2. If $u$ is a nonnegative solution of ${\cal L}u=f\ge0$ in 
${\cal B}_{\rho,\rho}$ satisfying $u\ge k$ on $\partial {\cal B}_{\rho,\rho}\cap\{x_n=\rho\}$ for some $k>0$ 
then for $\xi\le\frac 12$ the inequality (\ref{lem22}) holds in 
${\cal B}_{(1-\xi)\rho,\rho}\setminus{\cal B}_{(1-\xi)\rho,\xi\rho}$.
}\medskip

{\bf Proof}. This statement is more simple than Lemma 2.2. Consider the set 
$\widetilde K_{\rho}={\cal B}_{\rho,\rho}\setminus{\cal B}_{\rho,\frac {\xi\rho}{2}}$. Since 
coefficients $b^{(2)}_i$ are bounded on this set and $\|{\bf b}^{(2)}\|_{n,\widetilde K_{\rho}}$ is
under control, we can apply standard boundary growth lemma, and the statement follows.\hfill$\square$\medskip

{\bf Remark 2}. If we replace the assumption $f\ge0$ by $f_-\in L_n({\cal B}_{\rho,\rho})$,
the estimate (\ref{lem22}) holds true with additional term 
$-N_2(n,\nu,{\mathfrak B},\omega,\xi)\rho\cdot
\|f_-\|_{n,{\cal B}_{\rho,\rho}}$
in the right-hand side. The proof runs without changes.\medskip

Now we prove the main result of this Section.\medskip

{\bf Theorem 2.3}
. {\it Let $\cal L$ be an operator 
of the form (\,{\bf NDE}) in ${\cal B}_{R,R}$, and let the conditions (\ref{ell}), (\ref{decomp1}) and (\ref{decomp2}) be satisfied. Suppose also that for $\rho\le R$
\begin{equation}\label{Dini}
\|b^{(1)}_n\|_{n,{\cal B}_{\rho,\rho}}\le {\mathfrak B}_1\omega(\rho/R).
\end{equation}
Then\medskip

{\bf 1}. Any solution of ${\cal L}u=f\le0$ in ${\cal B}_{R,R}$ such that $u|_{x_n=0}\le0$
and $u(0)=0$ satisfies
$$\sup\limits_{0<x_n<R/2}\frac{u(0,x_n)}{x_n}\le \frac {N_3^+}R\cdot\sup\limits_{{\cal B}_{R/2,R/2}}u,$$ Consequently, if $D_nu(0)$ exists then $(D_nu)_+(0)$ is finite.\medskip

{\bf 2}. Any positive solution of ${\cal L}u=f\ge0$ in ${\cal B}_{R,R}$ such that $u(0)=0$
satisfies $\inf\limits_{0<x_n<R/2}\frac{u(0,x_n)}{x_n}>0$. Consequently, if $D_nu(0)$ exists, it is 
positive. If, in addition, $f\equiv0$, the following estimate holds:
$$\inf\limits_{0<x_n<R/2}\frac{u(0,x_n)}{x_n}\ge {N_3^-}\cdot \frac{u(0,R/2)}{R/2}.$$

The constants $N_3^{\pm}$ depend on $n$, $\nu$, ${\mathfrak B}$, ${\mathfrak B}_1$, $\omega$ and the the moduli of
continuity of $|{\bf b}|$ in $L_n({\cal B}_{R,R})$.}\medskip

{\bf Proof}. We introduce the sequence of cylinders ${\cal B}_{\rho_k,h_k}$, $k\ge0$, where $\rho_k=2^{-k}\rho_0$,
$h_k=\zeta_k\rho_k$, while $\rho_0\le R$ and the sequence $\zeta_k\downarrow0$ will be chosen later.

Denote by $M^{\pm}_k$, $k\ge1$, the quantities
$$M^+_k= \sup\limits_{{\cal B}_{\rho_k,h_{k-1}}}\frac{u(x)}{\max\{x_n,h_k\}}\ge
\sup\limits_{{\cal B}_{\rho_k,h_{k-1}}\setminus{\cal B}_{\rho_k,h_k}}\frac{u(x)}{x_n};\qquad
M^-_k=\inf\limits_{{\cal B}_{\rho_k,h_{k-1}}\setminus{\cal B}_{\rho_k,h_k}}\frac{u(x)}{x_n}.
$$
Note that in the case {\bf 2} $M^-_k>0$.

We define two function sequences
$$v^1_k=u-M^+_kh_k\,\frac {\varphi^+(x_n/R)}{\varphi^+(h_k/R)};\qquad v^2_k=M^-_kh_k\,\frac {\varphi^-(x_n/R)}{\varphi^-(h_k/R)}-u,$$
where, similarly to Lemma 2.2,
\begin{equation}\label{barrier1}
\varphi^{\pm}(s)=\int\limits_0^s\exp\Big(\mp\frac{\mathfrak B}{\nu}\,{\cal I}_{\omega}(\sigma)\Big) d\sigma,
\end{equation}
and denote 
$$\gathered
V_k=v^1_k;\quad M_k=M^+_k;\quad \Phi=\varphi^+ \qquad \mbox{in the case {\bf 1}};\\ 
V_k=v^2_k;\quad M_k=M^-_k;\quad \Phi=\varphi^- \qquad \mbox{in the case {\bf 2}}.
\endgathered$$
It is easy to see that $V_k\big|_{x_n=0}\le0$ while the definition of $M_k$ gives $V_k\le0$ 
on the top of the cylinder ${\cal B}_{\rho_k,h_k}$.

To estimate $V_k$, we refine a trick from \cite{S10}. Let $x^0\in {\cal B}_{\rho_k-h_k, h_k}$. 
Assume first that $x^0_n\le \frac{h_k}{2}$. Then we apply Lemma 2.2\footnote{To proceed we suppose 
that $\rho_0/R$ is so small that $\|{\bf b}\|_{n,{\cal B}_{\rho_0,\rho_0}}\le{\mathfrak B}_0$.} 
to the (positive) function $M_kh_k-V_k$ in ${\cal B}_{h_k, h_k}(x^0\vphantom{x}')$ (with regard to 
Remark 2). This gives for 
$x\in {\cal B}_{\frac {h_k}2,\frac {h_k}2}(x^0\vphantom{x}')$
\begin{multline}\label{qqq}
M_kh_k-V_k(x)\ge M_kh_k\cdot\big[\beta(n,\nu,{\mathfrak B},\omega,1/2)-\\
-N_1(n,\nu,{\mathfrak B},\omega,1/2)\cdot
(\|{\bf b}^{(1)}\|_{n,{\cal B}_{\rho_k,\rho_k}}+{\mathfrak B}\omega(\rho_k/R))\big]-\\
-N_2(n,\nu,{\mathfrak B},1/2)h_k\cdot\|({\cal L}V_k)_+\|_{n,{\cal B}_{h_k,h_k}(x^0\vphantom{x}')}.
\end{multline}
We suppose that $\rho_0/R$ is so small that the quantity in the square brackets is greater that $\frac {\beta}2$.
Further, direct calculation similar to Lemma 2.2 shows that the assumptions of theorem imply
$${\cal L}V_k\le M_k|b^{(1)}_n|\,\Phi'(x_n/R)\,\frac {h_k/R}{\Phi(h_k/R)}\qquad\mbox{in}
\quad {\cal B}_{\rho_k,h_k}.$$
Note that $\varphi^+$ is concave, $\varphi^-$ is convex, and both of them are increasing. Therefore,
$$\Phi'(x_n/R)\,\frac {h_k/R}{\Phi(h_k/R)}\le\frac {\max\{1,\Phi'(1)\}}{\Phi(1)}.$$
Substituting these inequalities into (\ref{qqq}) we arrive at
$$V_k(x)\le M_kh_k\cdot\Big[1-\beta/2+C_8(n,\nu,{\mathfrak B},\omega)
\|b^{(1)}_n\|_{n,{\cal B}_{h_k,h_k}(x^0\vphantom{x}')}\Big]\qquad\mbox{for}\quad
x\in{\cal B}_{\frac {h_k}2,\frac {h_k}2}(x^0\vphantom{x}').
$$
In particular, this estimate is valid for $x=x^0$. If $x^0_n\ge \frac{h_k}{2}$, we get the same estimate using Lemma 2.2' instead of Lemma 2.2.

Taking supremum w.r.t $x^0$, we obtain
$$\sup\limits_{{\cal B}_{\rho_k-h_k, h_k}}V_k\le M_kh_k\cdot\big(1-\beta/2+C_8{\mathfrak B}_1\omega(\rho_k/R)\big).
$$
Repeating previous arguments provides for $m\le \frac{\rho_k}{h_k}$
$$\sup\limits_{{\cal B}_{\rho_k-mh_k, h_k}}V_k\le
M_kh_k\cdot\Big((1-\beta/2)^m+C_8{\mathfrak B}_1\frac{\omega(\rho_k/R)}{\beta/2}\Big).
$$
Setting $m=\lfloor\frac{\rho_{k+1}}{h_k}\rfloor$, we arrive at
$$\sup\limits_{{\cal B}_{\rho_{k+1}, h_k}}V_k\le 
\frac{M_kh_k}{1-\beta/2}\cdot\Big(\exp\Big(-\lambda \frac{\rho_{k+1}}{h_k}\Big)+
C_8{\mathfrak B}_1\frac {\omega(\rho_k/R)}{\beta/2}\Big),
$$
where $\lambda=-\ln(1-\beta/2)>0$.

Therefore, for $x\in {\cal B}_{\rho_{k+1}, h_k}$
\begin{equation}\label{qqqq}
\frac{V_k(x)}{\max\{x_n,h_{k+1}\}}\le M_k\gamma_k,
\end{equation}
where $\gamma_k=\frac 1{1-\beta/2}\,\frac{\zeta_k}{2\zeta_{k+1}}\cdot
\big(\exp\big(- \frac{\lambda}{2\zeta_k}\big)+C_8{\mathfrak B}_1\frac {\omega(\rho_k/R)}{\beta/2}\big)$.

Estimate (\ref{qqqq}) implies in the cases {\bf 1} and {\bf 2}, respectively,
\begin{equation}\label{qqqqq}
M^{\pm}_{k+1}\lessgtr M^{\pm}_k(\delta^{\pm}_k\pm\gamma_k),
\end{equation}
where
\begin{equation}\label{delta}
\gathered
\delta^+_k=\frac {h_k}{\varphi^+(h_k/R)}\cdot\sup\limits_{0\le x_n\le h_k}
\frac {\varphi^+(x_n/R)}{\max\{x_n,h_{k+1}\}}=\frac {h_k}{\varphi^+(h_k/R)}\cdot
\frac {\varphi^+(h_{k+1}/R)}{h_{k+1}};\\
\delta^-_k=\frac {h_k}{\varphi^-(h_k/R)}\cdot\inf\limits_{h_{k+1}\le x_n\le h_k}
\frac {\varphi^-(x_n/R)}{x_n}=\frac {h_k}{\varphi^-(h_k/R)}\cdot\frac {\varphi^-(h_{k+1}/R)}{h_{k+1}}.
\endgathered
\end{equation}
Since $\lim\limits_{s\to0+}\frac {\Phi(s)}{s}=1$, we have
$$\prod\limits_k \delta^{\pm}_k=\frac {h_1/R}{\Phi(h_1/R)}\lessgtr\frac {1}{\Phi(1)}.
$$
Thus, (\ref{qqqqq}) gives in the cases {\bf 1} and {\bf 2}, respectively,
$$M^{\pm}_{k+1}\lessgtr \frac {M^{\pm}_1}{\Phi(1)}\cdot
\prod\limits_{j=1}^k\big(1\pm\frac {\gamma_j}{\delta^{\pm}_j}\big)\lessgtr \frac {M^{\pm}_1}{\Phi(1)}\cdot
\prod\limits_{j=1}^k(1\pm\gamma_j\cdot\max\{1,\Phi(1)\}).
$$
We set $\zeta_k=\frac 1{k+k_0}$ and choose $k_0$ so large and $\rho_0/R$ so small that 
$\gamma_1\cdot\max\{1,\Phi(1)\}\le \frac 12$. Note that $k_0$ and $\rho_0/R$ satisfying all the 
conditions imposed depend only on $n$, $\nu$, $\omega$, ${\mathfrak B}$ and ${\mathfrak B}_1$.

Now we observe that the first terms in $\gamma_k$ form a convergent series. The same is true for the second terms, since
$$\sum\limits_{k=1}^{\infty}\omega\big(2^{-k}\rho_0/R\big)
\asymp\int\limits_0^{\infty}\omega\big(2^{-s}\rho_0/R\big)\,ds\asymp{\cal I}_{\omega}(\rho_0/R).
$$
Thus, the series $\sum_k\gamma_k$ converges. Therefore, the infinite products 
$\Pi^{\pm}=\prod\limits_k(1\pm\gamma_k\cdot\max\{1,\Phi(1)\})$ also converge, and we obtain 
in the cases {\bf 1} and {\bf 2}, respectively,
$$M^{\pm}_{k}\lessgtr \frac {\Pi^{\pm}M^{\pm}_1}{\Phi(1)},\qquad k>1.
$$
Thus, all $M^+_k$ are bounded in the case {\bf 1}, and all $M^-_k$ are separated from zero in the case 
{\bf 2}.\medskip

Further, we note that $M^+_1\le \frac 1{h_1}\sup\limits_{{\cal B}_{R/2,R/2}}u$. This completes the proof
of the statement {\bf 1}.

If $f\equiv0$, then we set $K={\cal B}_{R,R}\setminus{\cal B}_{R,h_1/2}$. Similarly to Lemma 2.2, 
$\|{\bf b}^{(1)}+{\bf b}^{(2)}\|_{n,K}$ is bounded. Therefore, by the Harnack inequality 
(\cite[Theorem 3.3]{S10}), $M^{-}_1\asymp \frac{u(0,R/2)}{R/2}$. This completes the proof of the statement {\bf 2}.\hfill$\square$\medskip

{\bf Remark 3}. If we replace in the case {\bf 1} the assumption $f\le0$ by $f=f^{(1)}+f^{(2)}$ with 
$$\|f^{(1)}_+\|_{n,{\cal B}_{\rho,\rho}}\le {\mathfrak F}_1\omega(\rho/R), \qquad
f^{(2)}_+\le {\mathfrak F}_2\,\frac {\omega(x_n/R)}{x_n},
$$
the estimate 
$$\sup\limits_{0<x_n<R/2}\frac{u(0,x_n)}{x_n}\le 
N_3^+\cdot\Big(\frac 1R\,\sup\limits_{{\cal B}_{R/2,R/2}}u+{\mathfrak F}_1+{\mathfrak F}_2\Big)$$
remains valid. The proof runs with minor changes.\medskip

Let us compare Theorem 2.3 with results known earlier. Surely, the proof of ({\bf ZHO}) for classical
solutions works also for strong solutions if we apply the Aleksandrov maximum principle 
(\cite{Al}; see also a survey \cite{N05}, where the history of this topic is presented). 
So, it was known long ago for $b_i\in L_{\infty}$. 

In \cite{Li85} the Hopf--Oleinik Lemma was proved for classical solutions of ({\bf NDE}) in 
${\cal C}^{1+{\cal D}}$ domains. Note that in this case one can locally rectify $\partial\Omega$ 
using the regularized distance (\cite[Theorem 2.1]{Li85}). After this Theorem 4.1 \cite{Li85} follows 
from a particular case ${\bf b}^{(1)}\equiv0$ of Theorem 2.3, part {\bf 2}. Similarly, the boundary 
gradient estimates obtained in \cite{Li86} can be reduced to the same particular case of Theorem 2.3, part {\bf 1}.

The boundary gradient estimates for solution to ({\bf NDE}) were established in \cite{LU88} provided 
${\bf b}\in L_q$, $q>n$; the Hopf--Oleinik Lemma under the same condition was announced in \cite{NU}. 
In \cite{S10} the second part of Theorem 2.3 is proved for ({\bf NDE}) under assumption ${\bf b}\in L_n$, $b_n\in L_q$, $q>n$. 
In \cite{AN} the first part of Theorem 2.3 was proved for composite coefficients with $\omega(\sigma)=\sigma^{\alpha}$,
$\alpha\in\,]0,1[$.

To compare our result with \cite{S08}, we need an auxiliary statement.\medskip

{\bf Lemma 2.4}. {\it Let $\Psi:\,[0,\sigma_0]\to{\mathbb R}_+$ be a nondecreasing function.
Then there exist nondecreasing ${\cal C}^1$ functions $\Psi^{\pm}:\,[0,\sigma_0]\to{\mathbb R}_+$ 
such that $\Psi^-\le\Psi\le\Psi^+$, and

1) if $\Psi(\sigma)/\sigma^2$ is summable then $(\Psi^+)'(\sigma)/\sigma$ is summable;

2) if $\Psi(\sigma)/\sigma^2$ is nonsummable then $(\Psi^-)'(\sigma)/\sigma$ is nonsummable.}\medskip

{\bf Proof}. Without loss of generality we can assume $\sigma_0=1$. Consider the function 
$\Psi_1(\tau)=\Psi(\tau^{-1})$ and note that $\Psi_1$ is summable on $[1,+\infty[$ iff
$\Psi(\sigma)/\sigma^2$ is summable on $]0,1]$.

Now we define $\Psi_2(\tau)=\Psi_1(\lfloor\tau\rfloor)\cdot(\lceil\tau\rceil-\tau)+
\Psi_1(\lceil\tau\rceil)\cdot(\tau-\lfloor\tau\rfloor)$. Using the Cauchy convergence criterion, 
it is easy to check that $\tau\Psi_2'(\tau)$ is summable iff $\Psi_1$ is summable. 
Also it is evident that $\Psi_2(\tau+2)<\Psi_1(\tau)<\Psi_2(\tau-2)$.

Finally, we mollify $\Psi_2$ so that $\widetilde\Psi_2(\tau+2)<\Psi_1(\tau)<\widetilde\Psi_2(\tau-2)$
and set $\Psi^{\pm}(\sigma)=\widetilde\Psi_2(\sigma^{-1}\mp2)$, expanding $\Psi^+$ to $]\frac 13,1]$
in a proper way.\hfill$\square$\medskip

In \cite[Theorem 1.8]{S08} the Hopf--Oleinik Lemma was proved for solution to ({\bf NDE}) with 
${\bf b}\equiv0$ under assumption $0\in\partial\Omega$ and $\Omega\supset Q_{\Psi}$, where
$$Q_{\Psi}=\{x\in{\mathbb R}^n:\ |x'|\le \sigma_0,\ \Psi(|x'|)<x_n<\sigma_0\},
$$
while $\Psi(\sigma)/\sigma^2$ is summable.

By Lemma 2.4, this case can be reduced to $\Omega=Q_{\Psi^+}$. Then we again rectify the boundary and
use part 2 of Theorem 2.3. In the same manner, \cite[Theorem 1.9]{S08} follows from part 1 of Theorem 
2.3.\medskip

{\bf Remark 4}. Note that the assumption (\ref{Dini}) cannot be removed. Let us describe
corresponding counterexample (see also \cite{NU} and \cite{S10}).

Let $u(x)=x_n\cdot \ln^\alpha (|x|^{-1})$ in ${\cal B}_{R,R}$. Then direct calculation shows that $u$ 
satisfies an equation
$$-\Delta u+b_n(x)D_nu=0\qquad\mbox{with} \quad 
|b_n|\le \frac{C(\alpha)}{|x|\ln(|x|^{-1})}\in L_n({\cal B}_{R,R}),$$
if $R$ is small enough, and $u>0=u|_{x_n=0}$ in ${\cal B}_{R,R}$. However, it is easy to see that $D_nu(0)=0$ for $\alpha<0$ and $D_nu(0)=+\infty$ for $\alpha>0$.\medskip

The condition (\ref{decomp2}) is also sharp. A simple one-dimensional counterexample is given in
 \cite{M11}: the function 
$\phi(s)=\int\limits_0^s\exp\Big(-\int\limits_t^1\frac {\omega(\tau)}{\tau}\,d\tau\Big) dt$ 
(cf. (\ref{barrier1})) is positive on $]0,1]$, vanishes at zero and satisfies the equation
$$-\phi''(s)+\frac {\omega(s)}{s}\,\phi'(s)=0.
$$
However, if $\omega$ is not Dini continuous at zero then $\phi'(0)=0$.\medskip

We describe also a more rich family of counterexamples generalizing \cite[Theorem 1.11]{S08}.\medskip

Let $\Omega$ be a convex domain. Suppose that $\partial\Omega=\{x:\ x_n=\phi(x')\}$ at a neighborhood of 
the origin, $\phi\in {\cal C}^1$, $D'\phi(0)=0$, and $\omega(\rho)=\sup\limits_{|x'|\le\rho}|D'\phi(x')|$
is not Dini continuous at zero. Let $\cal L$ be an operator of the form ({\bf NDE}) in $\Omega$ with 
${\bf b}=0$, and let the condition (\ref{ell}) be satisfied. 

It is shown in \cite{AN11} that {\it any} solution of ${\cal L}u=0$ positive in $\Omega$ and vanishing 
on $\partial\Omega$ at a neighborhood of the origin satisfies $\sup\limits_{|x|<\rho}\frac{u(x)}{\rho}\to0$
as $\rho\to0$.

Now we rectify $\partial\Omega$ at a neighborhood of the origin using the regularized distance 
and obtain a uniformly elliptic operator of the form ({\bf NDE}) in ${\cal B}_{R,R}$ with 
$|{\bf b}(x)|\le {\mathfrak B}\,\frac {\omega(x_n/R)}{x_n}$ for which the Hopf--Oleinik lemma fails.

\section{Parabolic case}

In this section we assume $\partial_tu$, $D(Du)\in \widehat L_{q,\ell,loc}(Q)$ with some $q,\ell<\infty$ 
such that $\frac nq+\frac 1{\ell}=1$.

We recall an estimate which is a particular case of the statement in \cite[Sec.3]{N87}. For the isotropic
case it was proved in \cite{Kr86}.\medskip

{\bf Proposition 3.1}. {\it Let $\cal M$ be an operator of the form (\,{\bf NDP}) in a cylinder 
$Q\subset B_R\times\,]0,T[$, and let the condition (\ref{ell}) be satisfied. Suppose that 
${\bf b}\in L_{\infty}(Q)$, and a function $\mathbb B$ such that $\partial_t{\mathbb B}$, 
$D(D\mathbb B)\in L_{\infty}(Q)$ satisfies ${\cal M}{\mathbb B}\ge |{\bf b}|$ a.e. in $Q$. Then for
any solution of ${\cal M}u=f$ in $Q$ satisfying $u\big|_{\partial'Q}\le 0$, 
the following estimate holds:}
$$u\le N(n)\cdot\Big(\frac{\|{\mathbb B}\|_{\infty,Q}+R}{\nu}\Big)^{\frac nq}\cdot 
|\!|\!|f_+|\!|\!|_{q,\ell,\{u>0\}}.
$$

The next statement generalizes \cite[Theorem 2]{AN}. For the isotropic case it was proved in 
\cite[Theorem 5.2]{Li00}.\medskip

{\bf Lemma 3.1}. {\it Let $\cal M$ be an operator of the form (\,{\bf NDP}) in a cylinder ${\cal Q}_{R,R}$, 
$R\le1$, and let the condition (\ref{ell}) be satisfied. Suppose also that the vector function ${\bf b}$ can be 
written as follows:
\begin{equation}\label{decomp1'}
{\bf b}={\bf b}^{(1)}+{\bf b}^{(2)};\qquad |{\bf b}^{(1)}|\in \widehat L_{q,\ell}({\cal Q}_{R,R}),\quad 
\frac nq+\frac 1{\ell}=1,\quad q,\ell<\infty,
\end{equation}
\begin{equation}\label{decomp2'}
|{\bf b}^{(2)}(x;t)|\le {\mathfrak B}\,\frac {\omega(x_n/R)}{x_n},\qquad  \omega\in {\cal D}_1.
\end{equation}

Then for any solution of ${\cal M}u=f$ in ${\cal Q}_{R,R}$ satisfying $u\big|_{\partial'{\cal Q}_{R,R}}\le 0$, 
the following estimate holds:
\begin{equation}\label{estmax}
u\le N_4\cdot\big(|\!|\!|{\bf b}^{(1)}|\!|\!|^{\ell}_{q,\ell,{\cal Q}_{R,R}}+R\big)^{\frac nq}\cdot
|\!|\!|f_+|\!|\!|_{q,\ell,\{u>0\}},
\end{equation}
where $N_4$ depends only on $n$, $\nu$, $\ell$, ${\mathfrak B}$ and $\omega$.}

{\bf Proof}. We consider a sequence of operators
$${\cal M}_{\ep}\equiv \partial_t-a_{ij\ep}(x;t)D_iD_j+[b^{(1)}_{i\ep}(x;t)+b^{(2)}_{i\ep}(x;t)]D_i.
$$
Here $a_{ij\ep}$ are smooth functions satisfying (\ref{ell}) uniformly w.r.t. $\ep$ and tending to $a_{ij}$ a.e.
in ${\cal Q}_{R,R}$ as $\ep\to0$. Further,
$$b^{(1)}_{i\ep}(x;t)=\max\{|b^{(1)}_i(x;t)|; \ep^{-1}\}\cdot{\rm sign}(b^{(1)}_i(x;t));
$$
$$b^{(2)}_{i\ep}(x;t)=\begin{cases}b^{(2)}_i(x;t), & x_n>\ep;\\
b^{(2)}_i(x',\ep;t), & x_n\le\ep.
\end{cases}
$$

Now we consider the boundary value problem
$${\cal M}_{\ep}v=\Big(\frac {2{\mathfrak B}}{\nu}\ {\cal I}_{\omega}(1)+1\Big)\,|{\bf b}^{(1)}_{\ep}|+
{\mathfrak B}\ \frac {\omega(\rho/R)}{\rho}\cdot\Big(\frac {2{\mathfrak B}}{\nu}\ {\cal I}_{\omega}(1)-1\Big)_+
\quad\mbox{in}\quad {\cal Q}_{R,R};\qquad v\big|_{\partial'{\cal Q}_{R,R}}=0,
$$
where $\rho\le R$ will be chosen later. Denote by ${\mathbb B}^{(1)}_{\ep}$ a unique solution of this BVP. By the 
maximum principle (\cite{Kr}), ${\mathbb B}^{(1)}_{\ep}\ge0$. Define
$${\mathbb B}_{\ep}(x;t)={\mathbb B}^{(1)}_{\ep}(x;t)+
\frac {2{\mathfrak B}}{\nu}\ R\int\limits_{x_n/R}^1{\cal I}_{\omega}(s)\,ds.
$$
Then
$${\cal M}_{\ep}{\mathbb B}_{\ep}={\cal M}_{\ep}{\mathbb B}^{(1)}_{\ep}+
2{\mathfrak B}\frac {a_{nn\ep}}{\nu}\,\frac {\omega(x_n/R)}{x_n}-
[b^{(1)}_{n\ep}+b^{(2)}_{n\ep}]\cdot\frac {2{\mathfrak B}}{\nu}\ {\cal I}_{\omega}(x_n/R)\ge 
|{\bf b}^{(1)}_{\ep}|+|{\bf b}^{(2)}_{\ep}|+F(x),
$$
where
$$F(x)={\mathfrak B}\,\frac {\omega(x_n/R)}{x_n}\cdot
\Big(1-\frac {2{\mathfrak B}}{\nu}\ {\cal I}_{\omega}(x_n/R)\Big)+
{\mathfrak B}\ \frac {\omega(\rho/R)}{\rho}\cdot
\Big(\frac {2{\mathfrak B}}{\nu}\ {\cal I}_{\omega}(1)-1\Big)_+.
$$

We set $\rho=\min\{1;\widehat s\}R$, where $\widehat s$ is the root of 
${\cal I}_{\omega}(s)=\frac {\nu}{2{\mathfrak B}}$. Then, for $0<x_n<\rho$, we have
$$F(x)\ge {\mathfrak B}\,\frac {\omega(x_n/R)}{x_n}\cdot
\Big(1-\frac {2{\mathfrak B}}{\nu}\ {\cal I}_{\omega}(\rho/R)\Big)\ge0.
$$
If, otherwise, $\rho\le x_n\le R$, then
$$F(x)\ge{\mathfrak B}\,\frac {\omega(\rho/R)}{\rho}\cdot
\Big(1-\frac {2{\mathfrak B}}{\nu}\ {\cal I}_{\omega}(1)\Big)+
{\mathfrak B}\ \frac {\omega(\rho/R)}{\rho}\cdot\Big(\frac {2{\mathfrak B}}{\nu}\ {\cal I}_{\omega}(1)-1\Big)_+\ge0.
$$
So, in any case ${\cal M}_{\ep}{\mathbb B}_{\ep}\ge|{\bf b}^{(1)}_{\ep}|+|{\bf b}^{(2)}_{\ep}|$.

Using Proposition 3.1, we obtain the estimate
\begin{equation}\label{Bep}
u\le N(n)\cdot\Big(\frac{\|{\mathbb B}_{\ep}\|_{\infty,Q}+R}{\nu}\Big)^{\frac nq}\cdot 
|\!|\!|({\cal M}_{\ep}u)_+|\!|\!|_{q,\ell,\{u>0\}}
\end{equation}
for any function $u$ satisfying the assumptions of Lemma.

Obviously, 
\begin{equation}\label{Bep1}
\|{\mathbb B}_{\ep}\|_{\infty,Q}\le\|{\mathbb B}^{(1)}_{\ep}\|_{\infty,Q}+\frac {2{\mathfrak B}}{\nu}\ R{\cal I}_{\omega}(1).
\end{equation}
Further, the function ${\mathbb B}^{(1)}_{\ep}$ itself satisfies the assumptions of Lemma. Therefore, one can
set $u={\mathbb B}^{(1)}_{\ep}$ in (\ref{Bep}) and use (\ref{Bep1}) arriving at
\begin{equation}\label{B1ep}
\|{\mathbb B}^{(1)}_{\ep}\|_{\infty,Q}\le 
N(n)\cdot\Big(\frac{\|{\mathbb B}^{(1)}_{\ep}\|_{\infty,Q}+R(1+\frac {2{\mathfrak B}}{\nu}\ {\cal I}_{\omega}(1))}{\nu}\Big)^{\frac nq}
\cdot |\!|\!|({\cal M}_{\ep}{\mathbb B}^{(1)}_{\ep})_+|\!|\!|_{q,\ell,{\cal Q}_{R,R}}
\end{equation}
(we recall that ${\mathbb B}^{(1)}_{\ep}\ge0$).

If $\|{\mathbb B}^{(1)}_{\ep}\|_{\infty,Q}>R(1+\frac {2{\mathfrak B}}{\nu}\ {\cal I}_{\omega}(1))$ then (\ref{B1ep}) gives
$$\frac{\|{\mathbb B}^{(1)}_{\ep}\|_{\infty,Q}+R(1+\frac {2{\mathfrak B}}{\nu}\ {\cal I}_{\omega}(1))}{\nu}\le
\Big(\frac {2N(n)}{\nu}\cdot |\!|\!|({\cal M}_{\ep}{\mathbb B}^{(1)}_{\ep})_+|\!|\!|_{q,\ell,{\cal Q}_{R,R}}\Big)
^{\ell}
$$
(here we used $\frac nq+\frac 1{\ell}=1$). Thus, in any case we have
$$\frac{\|{\mathbb B}^{(1)}_{\ep}\|_{\infty,Q}+R(1+\frac {2{\mathfrak B}}{\nu}\ {\cal I}_{\omega}(1))}{\nu}\le
\Big(\frac {2N(n)}{\nu}\cdot |\!|\!|({\cal M}_{\ep}{\mathbb B}^{(1)}_{\ep})_+|\!|\!|_{q,\ell,{\cal Q}_{R,R}}\Big)
^{\ell}+\frac{2R(1+\frac {2{\mathfrak B}}{\nu}\ {\cal I}_{\omega}(1))}{\nu}.
$$
Substituting this estimate into (\ref{Bep}) and taking into account the definition of ${\mathbb B}^{(1)}_{\ep}$
we obtain (\ref{estmax}) for ${\cal M}_{\ep}$ instead of ${\cal M}$. Passage to the limit as $\ep\to0$ completes
the proof.\hfill$\square$\medskip

The next Lemma is parabolic analog of Lemmas 2.2 and 2.2'.\medskip

{\bf Lemma 3.2}. {\it Let $\cal M$ be an operator of the form (\,{\bf NDP}) in 
${\cal Q}_{\rho,\rho}$, $\rho\le R$, and let the conditions (\ref{ell}), (\ref{decomp1'}) and (\ref{decomp2'})
be satisfied. Suppose in addition that
\begin{equation}\label{Morrey}
|\!|\!|{\bf b}^{(1)}|\!|\!|_{q,\ell,{\cal Q}_{\rho,\rho}}\le {\mathfrak A}\rho^{\frac 1{\ell}}.
\end{equation}
Let $u$ be a nonnegative solution of ${\cal M}u=f\ge0$ in ${\cal Q}_{\rho,\rho}$.\medskip

{\bf 1}. If $u\ge k$ on $\partial'{\cal Q}_{\rho,\rho}\cap\{x_n=0\}$ for some $k>0$ then for 
$\xi\le\frac 12$ the inequality
\begin{equation}\label{lem32}
u\ge k\cdot\Big(\widehat\beta(n,\nu,\ell,{\mathfrak B},{\mathfrak A},\omega,\xi)-
N_5(n,\nu,\ell,{\mathfrak B},{\mathfrak A},\omega,\xi)\cdot
\big(\rho^{-\frac 1{\ell}}|\!|\!|{\bf b}^{(1)}|\!|\!|_{q,\ell,{\cal Q}_{\rho,\rho}}
+{\mathfrak B}\omega(\rho/R)\big)\Big)
\end{equation}
holds in ${\cal Q}_{(1-\xi)\rho,(1-\xi)\rho}$ with some positive constants $\widehat\beta$ and $N_5$.\medskip

{\bf 2}. If $u\ge k$ on $\partial'{\cal Q}_{\rho,\rho}\cap\{x_n=\rho\}$ for some $k>0$ then for 
$\xi\le\frac 12$ the inequality (\ref{lem32}) holds in 
${\cal Q}_{(1-\xi)\rho,\rho}\setminus{\cal Q}_{(1-\xi)\rho,\xi\rho}$.
}\medskip

{\bf Proof}. We prove the first statement. The proof of the second one is more simple, and we omit
it.\medskip

First, let $\xi=\frac 12$. Consider the barrier function $\widehat w(x;t)=w(x)+\frac t{\rho^2}$, where $w$ is 
defined in (\ref{barrier}) with a constant $A\ge\frac {\sqrt{n-1}}{\nu}+\frac 1{4\sqrt{n-1}}$ to be determined
later.\medskip

Similarly to Lemma 2.2, direct calculation shows that
${\cal M}\widehat w\le C_9(A,\nu,{\mathfrak B})|{\bf b}^{(1)}|\rho^{-1}$ in 
${\cal Q}_{\rho,\frac{\rho}A}$.\medskip

Further, $\widehat w(x;t)\le0$ on $\partial'{\cal Q}_{\rho,\frac{\rho}A}\setminus \{x_n=0\}$.
Finally, $\widehat w(x;t)\le1$ on $\partial'{\cal Q}_{\rho,\frac{\rho}A}\cap \{x_n=0\}$. This gives 
$k\widehat w-u\le0$ on $\partial'{\cal Q}_{\rho,\frac{\rho}A}$.\medskip

Lemma 3.1, condition (\ref{Morrey}) and relation $\frac nq+\frac 1{\ell}=1$ give for 
$(x;t)\in{\cal Q}_{\rho,\frac{\rho}A}$
$$u(x;t)\ge k\widehat w(x;t)-kN_4
\cdot\big({\mathfrak A}^{\ell}+1\big)^{\frac nq}\rho^{\frac nq}
\cdot|\!|\!|({\cal M}\widehat w)_+|\!|\!|_{q,\ell,{\cal Q}_{\rho,\frac{\rho}A}}\ge
k\widehat w(x;t)-kC_{10}\rho^{-\frac 1{\ell}}|\!|\!|{\bf b}^{(1)}|\!|\!|_{q,\ell,{\cal Q}_{\rho,\rho}},
$$
where $N_4$ is the constant from Lemma 3.1 while $C_{10}$ depends only on $n$, $\nu$, $\ell$, $A$, ${\mathfrak B}$ 
and ${\mathfrak A}$.\medskip

Similarly to Lemma 2.2, one can choose 
$A=A(n,\nu,{\mathfrak B},\omega)\ge\frac {\sqrt{n-1}}{\nu}+\frac 1{4\sqrt{n-1}}$ so large that 
$2(1+A)\varphi(1/A)\le\frac 1{100}$. 
Then direct calculation gives
\begin{equation}\label{p}
u\ge k\cdot\Big(\frac 1{20}-C_{11}(n,\nu,\ell,{\mathfrak B},{\mathfrak A},\omega)
\rho^{-\frac 1{\ell}}|\!|\!|{\bf b}^{(1)}|\!|\!|_{q,\ell,{\cal Q}_{\rho,\rho}}\Big)_+
\equiv k_2\qquad\mbox{in}\quad {\cal B}_{\frac {\rho}2,\frac {\rho}{10A}}\times\,]-\rho^2/2,0[.
\end{equation}

Now we consider the set 
$\widehat K_{\rho}={\cal Q}_{\frac {3\rho}4,\rho}\setminus{\cal Q}_{\frac {3\rho}4,\frac {\rho}{20A}}$. Note that 
coefficients $b^{(2)}_i$ are bounded on this set, and
\begin{equation}\label{pp}
|\!|\!|({\bf b}^{(2)}|\!|\!|_{q,\ell,\widehat K_{\rho}}\le C_{12}(n){\mathfrak B}\rho^{\frac 1{\ell}}
\Big(\int\limits_{\frac 1{20A}}^1\Big(\frac {\omega(s\rho/R)}s\Big)^q ds\Big)^{\frac 1q}\le 
C_{13}(n,\nu,\ell,{\mathfrak B},\omega){\mathfrak B}\rho^{\frac 1{\ell}}\omega(\rho/R).
\end{equation}
We proceed as \cite[Lemma 3.2]{LU85} (where the isotropic case was considered) and obtain
$$u\ge k_2\cdot\Big(\widehat\varkappa(n,\nu,\ell,{\mathfrak B},{\mathfrak A},\omega)-
C_{14}(n,\nu,\ell,{\mathfrak B},{\mathfrak A},\omega)
\rho^{-\frac 1{\ell}}|\!|\!|{\bf b}^{(1)}+{\bf b}^{(2)}|\!|\!|_{q,\ell,\widehat K}\Big)
\qquad\mbox{in}\quad {\cal Q}_{\frac {\rho}2,\frac {\rho}2}\setminus{\cal Q}_{\frac {\rho}2,\frac {\rho}{10A}}.
$$
By (\ref{p}) and (\ref{pp}) the statement for $\xi=\frac 12$ follows.\medskip

For arbitrary $\xi<\frac 12$ we apply the obtained statement in cylinders 
${\cal Q}_{2\xi\rho,2\xi\rho}({x^0}\vphantom{x}';t^0)$ with $|{x^0}\vphantom{x}'|\le (1-2\xi)\rho$, 
$(1-4\xi^2)\rho^2\le t^0\le0$. We arrive at
$$u\ge k\cdot\Big(\frac {\widehat\varkappa}{20}
-C_{15}(n,\nu,\ell,{\mathfrak B},{\mathfrak A},\omega)
\big(\rho^{-\frac 1{\ell}}|\!|\!|{\bf b}^{(1)}|\!|\!|_{q,\ell,{\cal Q}_{\rho,\rho}}
+{\mathfrak B}\omega(\rho/R)\big)\Big)_+\equiv k_3
$$
in ${\cal B}_{(1-\xi)\rho,\xi\rho}\times\,]-(1-3\xi^2)\rho^2,0[$.

Finally, as in the first step, one can proceed as \cite[Lemma 3.2]{LU85} in the set
${\cal Q}_{\rho,\rho}\setminus{\cal Q}_{\rho,\frac {\xi\rho}2}$,
and (\ref{lem32}) follows.\hfill$\square$\medskip

{\bf Remark 5}. If we replace the assumption $f\ge0$ by $f_-\in L_{q,\ell}({\cal Q}_{\rho,\rho})$,
the estimate (\ref{lem32}) holds true with additional term 
$-N_6(n,\nu,\ell,{\mathfrak B},{\mathfrak A},\omega,\xi)\rho^{\frac nq}\cdot
|\!|\!|f_-|\!|\!|_{q,\ell,{\cal Q}_{\rho,\rho}}$ in the right-hand side. The proof runs without changes.\medskip

{\bf Theorem 3.3}. {\it Let $\cal M$ be an operator 
of the form (\,{\bf NDP}) in ${\cal Q}_{R,R}$, and let the condition (\ref{ell}), (\ref{decomp1'}) 
and (\ref{decomp2'}) be satisfied. Suppose also that 
\begin{equation}\label{Morrey1}
{\mathfrak A}_1(\rho)\equiv
\sup\limits_{Q_{\rho}(x^0;t^0)\subset{\cal Q}_{R,R}}
\rho^{-\frac 1{\ell}}|\!|\!|{\bf b}^{(1)}|\!|\!|_{q,\ell,Q_{\rho}(x^0;t^0)}\to0,
\qquad \rho\to0,
\end{equation}
and for $\rho\le R$
\begin{equation}\label{Dini1}
\sup\limits_{{\cal Q}_{\rho,\rho}({x^0}\vphantom{x}';t^0)\subset{\cal Q}_{R,R}}
\rho^{-\frac 1{\ell}}|\!|\!|b^{(1)}_n|\!|\!|_{q,\ell,{\cal Q}_{\rho,\rho}({x^0}\vphantom{x}';t^0)}
\le {\mathfrak B}_1\omega(\rho/R).
\end{equation}
Then\medskip

{\bf 1}. Any solution of ${\cal M}u=f\le0$ in ${\cal Q}_{R,R}$ such that $u|_{x_n=0}\le0$
and $u(0;0)=0$ satisfies  
$$\sup\limits_{0<x_n<R/2}\frac{u(0,x_n;0)}{x_n}\le 
\frac {N_7^+}R\cdot\sup\limits_{{\cal Q}_{R/2,R/2}}u.$$
Consequently, if $D_nu(0;0)$ exists then $(D_nu)_+(0;0)$ is finite.\medskip

{\bf 2}. Any positive solution of ${\cal M}u=f\ge0$ in ${\cal Q}_{R,R}$ such that $u(0;0)=0$
satisfies $\inf\limits_{0<x_n<R/2}\frac{u(0,x_n;0)}{x_n}>0$. Consequently, if $D_nu(0;0)$ exists, it is 
positive. If, in addition, $f\equiv0$, the following estimate holds:
$$\inf\limits_{0<x_n<R/2}\frac{u(0,x_n;0)}{x_n}\ge {N_7^-}\cdot \frac{u(0,R/2;-R^2/2)}{R/2}.$$

The constants $N_7^{\pm}$ depend on $n$, $\nu$, $\ell$, ${\mathfrak B}$,
${\mathfrak B}_1$, ${\mathfrak A}_1$ and $\omega$.}\medskip

{\bf Remark 6}. If $|{\bf b}^{(1)}|\in \widehat L_{q,\widetilde\ell}({\cal Q}_{R,R})$ such that 
$q,\widetilde\ell<\infty$ and $\frac nq+\frac 2{\widetilde\ell}=1$ then (\ref{Morrey1}) is obviously satisfied.\medskip

{\bf Proof}. Similarly to Theorem 2.3, we introduce the sequence of cylinders ${\cal Q}_{\rho_k,h_k}$, $k\ge0$, 
where $\rho_k=2^{-k}\rho_0$, $h_k=\zeta_k\rho_k$, while $\rho_0\le R$ and the sequence $\zeta_k\downarrow0$ 
will be chosen later.

Denote by $\widehat M^{\pm}_k$, $k\ge1$, the quantities
$$\widehat M^+_k=\sup\limits_{{\cal Q}_{\rho_k,h_{k-1}}}\frac{u(x;t)}{\max\{x_n,h_k\}}\ge
\sup\limits_{{\cal Q}_{\rho_k,h_{k-1}}\setminus{\cal Q}_{\rho_k,h_k}}\frac{u(x;t)}{x_n};\qquad
\widehat M^-_k=\inf\limits_{{\cal Q}_{\rho_k,h_{k-1}}\setminus{\cal Q}_{\rho_k,h_k}}\frac{u(x;t)}{x_n}.
$$
Note that in the case {\bf 2} $\widehat M^-_k>0$.

We define two function sequences
$$\widehat v^1_k=u-\widehat M^+_kh_k\,\frac {\varphi^+(x_n/R)}{\varphi^+(h_k/R)};\qquad 
\widehat v^2_k=\widehat M^-_kh_k\,\frac {\varphi^-(x_n/R)}{\varphi^-(h_k/R)}-u,$$
where functions $\varphi^{\pm}$ are introduced in (\ref{barrier1}), and denote 
$$\gathered
\widehat V_k=\widehat v^1_k;\quad \widehat M_k=\widehat M^+_k;\quad \Phi=\varphi^+ \qquad \mbox{in the case {\bf 1}};\\ 
\widehat V_k=\widehat v^2_k;\quad \widehat M_k=\widehat M^-_k;\quad \Phi=\varphi^- \qquad \mbox{in the case {\bf 2}}.
\endgathered$$
It is easy to see that $\widehat V_k\big|_{x_n=0}\le0$ while the definition of $\widehat M_k$ gives 
$\widehat V_k\big|_{x_n=h_k}\le0$.

To estimate $\widehat V_k$, we consider $(x^0;t^0)\in {\cal Q}_{\rho_k-h_k, h_k}$. 
Let $x^0_n\le \frac{h_k}{2}$. Then we apply the first part of Lemma 3.2 to the function 
$\widehat M_kh_k-\widehat V_k$ in ${\cal Q}_{h_k, h_k}(x^0\vphantom{x}';t^0)$ (with regard to Remark 5). 
This gives for $x\in {\cal Q}_{\frac {h_k}2,\frac {h_k}2}(x^0\vphantom{x}';t^0)$
\begin{multline*}\label{ppp}
\widehat M_kh_k-\widehat V_k(x)\ge \widehat M_kh_k\cdot
\big[\widehat\beta(n,\nu,\ell,{\mathfrak B},{\mathfrak A}_1(h_k),\omega,1/2)-\\
-N_5(n,\nu,\ell,{\mathfrak B},{\mathfrak A}_1(h_k),\omega,1/2)\cdot
({\mathfrak A}_1(h_k)+{\mathfrak B}\omega(\rho_k/R))\big]-\\
-N_6(n,\nu,\ell,{\mathfrak B},{\mathfrak A}_1(h_k),\omega,1/2)h_k^{\frac nq}
\cdot|\!|\!|({\cal M}\widehat V_k)_+|\!|\!|_{q,\ell,{\cal Q}_{h_k,h_k}(x^0\vphantom{x}';t^0)}.
\end{multline*}

By (\ref{Morrey1}), we can choose $\rho_0/R$ is so small that the quantity in the square brackets is greater 
that $\frac {\widehat\beta}2$. As in Theorem 2.3, for $(x;t)\in{\cal Q}_{\frac {h_k}2,\frac {h_k}2}(x^0\vphantom{x}';t^0)$
we arrive at
$$\widehat V_k(x;t)\le \widehat M_kh_k\cdot\Big[1-\widehat\beta/2+
C_{16}(n,\nu,\ell,{\mathfrak B},{\mathfrak A}_1(h_k),\omega)
h_k^{-\frac 1{\ell}}|\!|\!|b^{(1)}_n|\!|\!|_{q,\ell,{\cal Q}_{h_k,h_k}(x^0\vphantom{x}';t^0)}\Big].
$$
In particular, this estimate is valid for $(x;t)=(x^0;t^0)$. If $x^0_n\ge \frac{h_k}{2}$, we get the same 
estimate using the second part of Lemma 3.2.

Taking supremum w.r.t $(x^0;t^0)$, we obtain
$$\sup\limits_{{\cal Q}_{\rho_k-h_k, h_k}}\widehat V_k\le \widehat M_kh_k\cdot\big(1-\widehat\beta/2+
C_{16}{\mathfrak B}_1\omega(h_k/R)\big).
$$
Repeating previous arguments provides for $m\le \frac{\rho_k}{h_k}$
$$\sup\limits_{{\cal Q}_{\rho_k-mh_k, h_k}}\widehat V_k\le
\widehat M_kh_k\cdot\Big((1-\widehat\beta/2)^m+C_{16}{\mathfrak B}_1\frac{\omega(h_k/R)}{\widehat\beta/2}\Big).
$$
Setting $m=\lfloor\frac{\rho_{k+1}}{h_k}\rfloor$, we arrive at
$$\sup\limits_{{\cal Q}_{\rho_{k+1}, h_k}}\widehat V_k\le 
\frac{\widehat M_kh_k}{1-\widehat\beta/2}\cdot\Big(\exp\Big(-\widehat \lambda \frac{\rho_{k+1}}{h_k}\Big)+
C_{16}{\mathfrak B}_1\frac {\omega(h_k/R)}{\widehat\beta/2}\Big),
$$
where $\widehat \lambda=-\ln(1-\widehat\beta/2)>0$.

Therefore, for $(x;t)\in {\cal Q}_{\rho_{k+1}, h_k}$
\begin{equation}\label{pppp}
\frac{\widehat V_k(x;t)}{\max\{x_n,h_{k+1}\}}\le \widehat M_k\widehat \gamma_k,
\end{equation}
where $\widehat \gamma_k=\frac 1{1-\widehat\beta/2}\,\frac{\zeta_k}{2\zeta_{k+1}}\cdot
\big(\exp\big(- \frac{\widehat \lambda}{2\zeta_k}\big)+C_{16}{\mathfrak B}_1\frac {\omega(h_k/R)}{\widehat\beta/2}\big)$.

Estimate (\ref{pppp}) implies in the cases {\bf 1} and {\bf 2}, respectively,
$$\widehat M^{\pm}_{k+1}\lessgtr \widehat M^{\pm}_k(\delta^{\pm}_k\pm\widehat \gamma_k),
$$
where $\delta^{\pm}_k$ are defined in (\ref{delta}). Similarly to Theorem 2.3, we obtain
$$\widehat M^{\pm}_{k+1}\lessgtr \frac {\widehat M^{\pm}_1}{\Phi(1)}\cdot
\prod\limits_{j=1}^k(1\pm\widehat \gamma_j\cdot\max\{1,\Phi(1)\}).
$$
We set $\zeta_k=\frac 1{k+k_0}$ and choose $k_0$ so large and $\rho_0/R$ so small that 
$\widehat \gamma_1\cdot\max\{1,\Phi(1)\}\le \frac 12$. Note that $k_0$ and $\rho_0/R$ satisfying all 
the conditions imposed depend only on $n$, $\nu$, ${\mathfrak B}$, ${\mathfrak B}_1$, ${\mathfrak A}_1$ 
and $\omega$.

Now, as in Theorem 2.3, we observe that the series $\sum_k\widehat\gamma_k$ converges. Therefore, the infinite products 
$\widehat\Pi^{\pm}=\prod\limits_k(1\pm\widehat\gamma_k\cdot\max\{1,\Phi(1)\})$ also converge, and we obtain 
in the cases {\bf 1} and {\bf 2}, respectively,
$$\widehat M^{\pm}_{k}\lessgtr \frac {\widehat\Pi^{\pm}\widehat M^{\pm}_1}{\Phi(1)},\qquad k>1.
$$
Thus, all $\widehat M^+_k$ are bounded in the case {\bf 1}, and all $\widehat M^-_k$ are separated from zero in the case {\bf 2}.\medskip

Further, we note that $\widehat M^+_1\le \frac 1{h_1}\sup\limits_{{\cal Q}_{R/2,R/2}}u$.
This completes the proof of the statement {\bf 1}.

If $f\equiv0$ then we set 
$\widehat K={\cal Q}_{R,R}\setminus{\cal Q}_{R,h_1/2}$. Similarly to Lemma 3.2, 
$$\sup\limits_{Q_{\rho}(x^0,t^0)\subset\widehat K}\rho^{-\frac 1{\ell}}
|\!|\!|{\bf b}^{(1)}+{\bf b}^{(2)}|\!|\!|_{q,\ell,Q_{\rho}(x^0,t^0)}\to0,\qquad \rho\to0.$$ 
Therefore, we use the Harnack inequality which can be proved in a similar way as \cite[Theorem 3.3]{S10} 
(for bounded lower-order terms see \cite{KS}) and obtain 
$$M^-_1\ge \frac {C_{17}}{h_1}\cdot u(0,R/2;-R^2/2),$$
where $C_{17}$ depends on $n$, $\nu$, $\ell$, ${\mathfrak B}$,
 ${\mathfrak A}_1$ and $\omega$. This completes the proof of the statement {\bf 2}.\hfill$\square$\medskip

{\bf Remark 7}. If we replace in the case {\bf 1} the assumption $f\le0$ by $f=f^{(1)}+f^{(2)}$ with 
$$\sup\limits_{{\cal Q}_{\rho,\rho}({x^0}\vphantom{x}';t^0)\subset{\cal Q}_{R,R}}
\rho^{-\frac 1{\ell}}|\!|\!|f^{(1)}_+|\!|\!|_{q,\ell,{\cal Q}_{\rho,\rho}({x^0}\vphantom{x}';t^0)}
\le {\mathfrak F}_1\omega(\rho/R), 
\qquad f^{(2)}_+\le {\mathfrak F}_2\,\frac {\omega(x_n/R)}{x_n},
$$
the estimate 
$$\sup\limits_{0<x_n<R/2}\frac{u(0,x_n)}{x_n}\le 
N_7^+\cdot\Big(\frac 1R\,\sup\limits_{{\cal Q}_{R/2,R/2}}u+{\mathfrak F}_1+{\mathfrak F}_2\Big)$$
remains valid. The proof runs with minor changes.\medskip

Let us compare Theorem 3.3 with results known earlier. As in elliptic case, the proof of 
Hopf--Oleinik Lemma for classical solutions to parabolic equations with $b_i\in L_{\infty}$ works also 
for strong solutions by the Aleksandrov--Krylov maximum principle (\cite{Kr}; see also \cite{N05}). 

In \cite{KH} the Hopf--Oleinik Lemma was proved for classical solutions of ({\bf NDP}) in 
${\cal C}^{1+{\cal D},\frac 12+{\cal D}}$ domains. Using the parabolic regularized distance 
(\cite[Theorem 3.1]{Li85}) one can locally rectify the boundary and reduce the result of \cite{KH} 
to a particular case ${\bf b}^{(1)}\equiv0$ of Theorem 3.3.

The boundary gradient estimates for solutions to ({\bf NDP}) were established in \cite{LU88} provided 
$|{\bf b}|\in L_{q+2}$, $q>n$; the Hopf--Oleinik Lemma under condition 
$|{\bf b}|\in \widehat L_{q,\widetilde\ell}$, $\frac nq+\frac 2{\widetilde\ell}<1$, 
$q,\widetilde\ell<\infty$, was announced in \cite{NU}. In \cite{AN} the first part of Theorem 3.3 
was proved for composite coefficients with
${\bf b}^{(1)}\in L_{q+2}$, $q>n$, and $\omega(\sigma)=\sigma^{\alpha}$, $\alpha\in\,]0,1[$.\medskip

{\bf Remark 8}. Note that the assumption (\ref{Dini1}) cannot be removed. Let us describe
corresponding counterexample (see \cite{NU}).

Let $u(x;t)=x_n\cdot \ln^\alpha ((|x|^2-t)^{-1})$ in ${\cal Q}_{R,R}$. Then direct calculation shows 
that $u$ satisfies an equation
$$\partial_tu-\Delta u+b_n(x;t)D_nu=0\qquad\mbox{with} \quad 
|b_n|\le \frac{C(\alpha)}{(|x|^2-t)^{\frac 12}\ln((|x|^2-t)^{-1})}\in 
\widehat L_{q,\widetilde\ell}({\cal Q}_{R,R}),$$
for any $q,\widetilde\ell<\infty$ such that $\frac nq+\frac 2{\widetilde\ell}=1$, if $R$ is small enough.
By Remark 6, the assumption (\ref{Morrey1}) is satisfied. Moreover, $u>0=u|_{x_n=0}$ in ${\cal Q}_{R,R}$.
 However, it is easy to see that $D_nu(0;0)=0$ for $\alpha<0$ and 
$D_nu(0;0)=+\infty$ for $\alpha>0$.\medskip

The condition (\ref{decomp2'}) is also sharp. Indeed, considering
functions depending only on spatial variables we see that the counterexample 
at the end of Section 2 works also for the parabolic operator ({\bf NDP}). A more rich family of
counterexamples also can be extracted from \cite{AN11}.


\begin{thebibliography}{AN95}

\bibitem[Al]{Al}
A.D.~Aleksandrov, {\em Uniqueness conditions and bounds for the solution of the Dirichlet problem},
Vestnik Leningrad. Univ. Ser. Mat. Meh. Astronom. {\bf 18} (1963), N3, 5--29 (Russian).

\bibitem[A-Z]{M11}
R.~Alvarado, D.~Brigham, V.~Maz'ya, M.~Mitrea, E.~Ziad\'e, {\em On the regularity of domains 
satisfying a uniform hour-glass condition and a sharp version of the Hopf--Oleinik boundary point principle},
Probl. Mat. Anal., {\bf 57} (2011), 3--68 (Russian); English transl.: J. Math. Sci.,
{\bf 176} (2011), N3, 281--360.

\bibitem[AN95]{AN}
D.E.~Apushkinskaya, A.I.~Nazarov, {\em Boundary estimates for the first-order 
derivatives of a solution to a nondivergent parabolic equation with composite 
right-hand side and coefficients of lower-order derivatives}, Probl. 
Mat. Anal., {\bf 14} (1995), 3--27 (Russian); English transl.: J. Math. Sci.,
{\bf 77} (1995), N4, 3257--3276.

\bibitem[AN11]{AN11}
D.E.~Apushkinskaya, A.I.~Nazarov, {\em A counterexample to the Hopf--Oleinik lemma}, in preparation.

\bibitem[G]{G}
G.~Giraud, {\em Probl\`emes de valeurs \`a la fronti\`ere relatifs \`a certaines donn\'es discontinues},
Bull. de la Soc. Math. de France, {\bf 61} (1933), 1--54.

\bibitem[GT]{GT}
D.~Gilbarg, N.~Trudinger, {\em Elliptic Partial Differential Equations of Second Order}, 2nd ed.,
Springer, Berlin etc. (1983).


\bibitem[Hi]{Hi}
B.N.~Him\v{c}enko, {\em On the behavior of solutions of elliptic equations near the boundary of a domain 
of type $A^{(1)}$}, DAN SSSR {\bf 193} (1970), 304--305 (Russian); English transl.: Sov. Math. Dokl. 
{\bf 11} (1970), 943--944.

\bibitem[Ho]{Ho}
E.~Hopf, {\em A remark on linear elliptic differential equations of second order}, Proc. AMS, 
{\bf 3} (1952), 791--793.

\bibitem[K]{K}
L.I.~Kamynin, {\em A theorem on the interior derivative for a second-order uniformly parabolic equation},
DAN SSSR {\bf 299} (1988), N2, 280--283 (Russian); English transl.: Sov. Math. Dokl. {\bf 37} (1988), N2, 373--376.

\bibitem[KHi]{KH}
L.I.~Kamynin, B.N.~Him\v{c}enko, {\em The analogues of the Giraud theorem for a second order 
parabolic equation}, Sibirsk. Mat. Zh. {\bf 14} (1973), N1, 86--110 (Russian); English transl.:
Siberian Math. J. {\bf 14} (1973), 59--77.

\bibitem[Kr76]{Kr}
N.V.~Krylov, {\em Sequences of convex functions, and estimates of the maximum of the solution of 
a parabolic equation}, Sibirsk. Mat. Zh. {\bf 17} (1976), N2, 290--303 (Russian);
English transl.: Siberian Math. J. {\bf 17} (1976), N2, 226--236.

\bibitem[Kr86]{Kr86}
N.V.~Krylov, {\em On estimates for the maximum of solutions of a parabolic equation and estimates for
distribution of a semimartingal}, Mat. Sb. {\bf 130(172)} (1986), N2(6), 207--221 (Russian); English transl.: Math. USSR -- Sbornik, {\bf 58} (1987), N1, 207--221.

\bibitem[KrS]{KS} 
N.V.~Krylov, M.V.~Safonov, {\em A property of the solutions of parabolic equations with measurable coefficients}, Izv. AN SSSR Ser. Mat. {\bf 44} (1980), N1, 161--175 (Russian); English transl.: 
Math. USSR -- Izvestiya, {\bf 16} (1981), N1, 151--164 


\bibitem[LU85]{LU85}
O.A.~Ladyzhenskaya, N.N.~Uraltseva, {\em Estimates of the H\"older constant for functions satisfying 
a uniformly elliptic or a uniformly parabolic quasilinear inequality with unbounded coefficients}, ZNS LOMI, 
{\bf 147} (1985), 72--94 (Russian); English transl.: J. Sov. Math., {\bf 37} (1987), N1, 837--851.

\bibitem[LU88]{LU88}
O.A.~Ladyzhenskaya, N.N.~Uraltseva, {\em Estimates on the boundary of a domain for the first derivatives
of functions satisfying an elliptic or parabolic inequality}, Trudy MI AN SSSR, {\bf 179} (1988), 102--125
(Russian); English transl.: Proc. Steklov Inst. Math., {\bf 179} (1989), 109--135.

\bibitem[Li85]{Li85}
G.M.~Lieberman, {\em Regularized distance and its applications}, Pacific J. Math. {\bf 117} (1985), N2, 329--352. 

\bibitem[Li86]{Li86}
G.M.~Lieberman, {\em The Dirichlet problem for quasilinear elliptic equations with continuously
differentiable boundary data}, Commun. PDEs, {\bf 11} (1986), N2, 167--229.

\bibitem[Li00]{Li00}
G.M.~Lieberman, {\em The maximum principle for equations with composite coefficients},
Electron. J. Diff. Eqs, N38 (2000), 17 pp. 

\bibitem[Na]{Na}
N.S.~Nadirashvili, {\em On the question of the uniqueness of the solution of the second boundary value
problem for second-order elliptic equations}, Mat. Sb. {\bf 122(164)} (1983), N3, 341--359 (Russian); 
English transl.: Math. USSR -- Sbornik, {\bf 50} (1985), N2, 325--341.


\bibitem[N87]{N87}
A.I.~Nazarov, {\em Interpolation of linear spaces and estimates for the maximum of 
a solution for parabolic equations}, Partial differential equations, Akad. 
Nauk SSSR Sibirsk. Otdel., Inst. Mat., No\-vo\-si\-birsk, 1987, 50--72 (Russian).

\bibitem[N01]{N01}
A.I.~Nazarov, {\em Estimates for the maximum of solutions of elliptic and 
parabolic equations in terms of weighted norms of the right-hand side}, 
Alg. \& Anal. {\bf 13} (2001), N2, 151--164 (Russian); English transl.: 
St.Petersburg Math. J. {\bf 13} (2002), N2, 269--279.

\bibitem[N05]{N05}
A.I.~Nazarov, {\em The A.D.~Aleksandrov maximum principle}, Contemp. Math. Appl. {\bf 29} (2005), 
127--143 (Russian); English transl.: J. Math. Sci., {\bf 142} (2007), N3, 2154--2171.

\bibitem[NU]{NU}
A.I. Nazarov, N.N. Ural'tseva, {\em Qualitative properties of solutions 
to elliptic and parabolic equations with unbounded lower-order coefficients},
St.Petersburg Math. Soc. El. Prepr. Archive. N~2009-05. 6pp.

\bibitem[Ni]{Ni}
L.~Nirenberg, {\em A strong maximum principle for parabolic equations}, Comm. Pure Appl. Math. {\bf 6}
(1953), 167--177.

\bibitem[O]{O}
O.A.~Oleinik, {\em On properties of solutions of certain boundary problems for equations of
elliptic type}, Mat. Sb. (N.S.), {\bf 30} (1952), 695--702 (Russian).

\bibitem[S08]{S08}
M.V.~Safonov, {\em Boundary estimates for positive solutions to second order elliptic equations}, preprint, http://arxiv.org/abs/0810.0522 (2008), 20pp.

\bibitem[S10]{S10}
M.V.~Safonov, {\em Non-divergence elliptic equations of second order with unbounded drift},
AMS Transl., Ser. 2, {\bf 229} (2010), 211--232.

\bibitem[Z]{Z}
M.S.~Zaremba, {\em Sur un probl\`eme mixte relatif \`a l'\'equation de Laplace}, Bull. Intern. de 
l'Acad. Sci. de Cracovie, Ser. A, Sci. Math. (1910), 313--344.


\end{thebibliography}
\end{document}